\documentclass[12pt]{amsart}
\usepackage{amssymb}
\usepackage[headings]{fullpage}
\usepackage{amsfonts}
\usepackage{paralist}
\usepackage{xcolor}
\usepackage{amsmath}

\usepackage[colorlinks=true,linkcolor=blue,citecolor=blue,urlcolor=blue,
pagebackref=true]{hyperref} 

\makeatletter
\@addtoreset{equation}{section}
\def\theequation{\thesection.\@arabic \c@equation}

\def\theenumi{\@roman\c@enumi}

\makeatother

\newtheorem{theorem}[equation]{Theorem}
\newtheorem{lemma}[equation]{Lemma}
\newtheorem{proposition}[equation]{Proposition}

\newtheorem{conjecture}[equation]{Conjecture}

\newtheorem{claim*}{Claim}

\theoremstyle{definition}
\newtheorem{remark}[equation]{Remark}

\newtheorem{eg}[equation]{Example}
\newenvironment{example}[1][]{\begin{eg}[#1] \pushQED{\qed}}{\popQED
\end{eg}}
\newtheorem{definition}[equation]{Definition}

\newtheorem{notn}[equation]{Notation}
\newenvironment{notation}[1][]{\begin{notn}[#1]\pushQED{\qed}}{\popQED
\end{notn}}

\newcommand{\bbB}{\mathbb{B}}

\newcommand{\naturals}{\mathbb N}
\newcommand{\ints}{\mathbb Z}
\newcommand{\define}{\emph}

\DeclareMathOperator{\Hom}{Hom}

\DeclareMathOperator{\codim}{codim}
\DeclareMathOperator{\projdim}{pd}
\DeclareMathOperator{\Tor}{Tor}

\newsavebox{\upperboundtheorem}

\title[Hilbert coefficients of Gorenstein modules]%
{An upper bound for the first Hilbert coefficient of Gorenstein algebras
and modules}

\author[S.~El~Khoury]{Sabine El Khoury} 
\address{Department of Mathematics, American University of Beirut, Beirut,
Lebanon.}
\email{se24@aub.edu.lb} 

\author{Manoj Kummini}
\address{Chennai Mathematical Institute, Siruseri,
Tamilnadu 603103 India.}
\email{mkummini@cmi.ac.in}

\author[H.~Srinivasan]{Hema Srinivasan}
\address{Department of Mathematics, University of Missouri, Columbia,
Missouri, USA.}
\email{hema@math.missouri.edu}

\begin{document}

\begin{abstract} Let $R$ be a polynomial ring over a field 
and $M= \bigoplus_n M_n$ a finitely generated
graded $R$-module, minimally generated by homogeneous elements of degree
zero with a graded $R$-minimal free resolution  $\mathbf{F}$. A
Cohen-Macaulay module $M$ is Gorenstein when the graded resolution is
symmetric.  We give an upper bound for the first Hilbert coefficient, $e_1$
in terms of the shifts in the graded resolution of $M$.   When $M = R/I$, a
Gorenstein algebra, this bound agrees with the bound obtained in \cite{ES}
in Gorenstein algebras with quasi-pure resolution.   We conjecture a
similar bound for the higher coefficients.
\end{abstract}
\maketitle

\section{Introduction}
Let $\Bbbk$ be any field and  $R=\bigoplus R_i$ be a 
polynomial ring over $\Bbbk$, finitely generated in degree one. 
Let $M = \bigoplus_i M_i$ be a finitely generated graded $R$-module.
We write $H(M,i ) := \dim_{\Bbbk}M_i$ for the 
\define{Hilbert function} of $M$.
For $i \gg 0$, $H(M,i) = P_M(i)$ where $P_M(x)$ is a polynomial
called the \define{the Hilbert polynomial} of $M$. 
The degree of  $P_M(x)$  is $\dim M-1$. We can write
  $$P_M(x)= \displaystyle \sum_{i=0}^{d-1} (-1)^{i}e_{i}  {  x+d-1-i
  \choose x} \\
= \frac {e_0}{(d-1)!} x^{d-1}+\ldots +(-1)^{d-1}e_{d-1}\\
$$
where the coefficients $e_i$ are non-negative integers (depending on $M$), 
called the \define{Hilbert coefficients} of $M$. 
The first one, $e_0$, also denoted by $e$, is the multiplicity of $M$. 

We now introduce the notion of a Gorenstein $R$-module, generalizing the
notion of a Gorenstein quotient ring of $R$. Assume that $M$ is generated
minimally by elements of degree $0$. Let
$$\mathbf{F} = \mathbf{F}(M) :  
0\rightarrow \displaystyle\bigoplus_{j=t_s}^{T_s}
  R(-j)^{\beta_{sj} }\stackrel{\delta_s}{\rightarrow} \ldots\rightarrow
  \displaystyle\bigoplus_{j=t_i}^{T_i}R(-j)^{\beta_{ij}}
  \stackrel{\delta_i}{\rightarrow}  \ldots\rightarrow
  \displaystyle\bigoplus_{j=t_1}^{T_1}R(-j)^{\beta_{1j}}\stackrel{\delta_1}
  \rightarrow R^{\beta_0}  \rightarrow M $$
be a minimal graded free resolution of $M$ as an $R$-module.
We assume that $\beta_{i,t_i } \neq 0$ and 
$\beta_{i,T_i } \neq 0$ for every $i$.
The \define{codimension} of $M$, denoted $\codim M$, is the height of the
annihilator of $M$. Note that $s \geq \codim M$, with equality if and only
if $M$ is Cohen-Macaulay. We say that $M$ is \define{Gorenstein} if it is
Cohen-Macaulay and if $\mathbf{F}$ is self-dual, i.e., $\Hom_R(\mathbf{F},
R) \simeq R$ after appropriate shifts.  

The numbers $\beta_{ij }$ appearing in $\mathbf{F}$ are called the
\define{graded Betti numbers} of $M$; note that $\beta_{i,j }$ is the 
number of copies of $R(-j)$ that appear at homological degree $i$ in
any minimal graded free resolution of $M$, and that 
$\beta_{i,j} = \dim_\Bbbk \Tor_i^R(\Bbbk, M)_j$.
If $M$ is Gorenstein, then
$t_s = T_s$ and $\beta_{ij} = \beta_{s-i, N-j}$ for  $N= T_s$.  
 
 There has been a lot of work on bounding the multiplicity and the Hilbert
 coefficients  of  finitely generated modules, in terms of the minimal and
 maximal shifts of their minimal free resolutions. C. Huneke and Srinivasan
 (see [\cite{HS}, Conjecture 1]) and of J. Herzog and Srinivasan
(see [\cite{HS}, Conjecture 2]) proposed upper and lower bounds for the
multiplicity in terms of $t_i$ and $T_i$. These conjectures were proven by
using the Boij and S\"oderberg in full. In \cite{HZ09}, Herzog and Zheng
extend the result on the multiplicity to all coefficients in the Cohen
Macaulay case, and found upper and lower bounds for the $e_i$ in the sense
of \cite{HS}. In the Gorenstein case, when the resolution is quasi-pure, i.e.,
$t_i \geq T_{i-1}$ for all $i$, the duality of the resolution is used to
find sharper bounds for the multiplicity and all Hilbert coefficients, see
\cite{S98} and \cite{ES1} respectively. When the resolution is not
quasi-pure, the duality of the resolution can be picked up in the Betti
table from the Boij and S\"oderberg decomposition. This gives a
generalization for the upper bound of the multiplicity to all Gorenstein
algebras, see \cite{EKS} for instance. In this paper, we extend the upper
bound found in  \cite{EKS} to the first coefficient $e_1$, and prove the
following theorem:

\begin{theorem} Let $M$ be a finitely generated graded Gorenstein
$R$-module, minimally generated by homogeneous elements of degree zero.
Let $s = \codim M$ and $k = \lfloor \frac{s}{2} \rfloor$.
Let $\beta_0(M)$ denote the minimal number of generators of $M$ and
$\mathbf{F}(M)$ be the graded resolution of $M$ as above.   Set  
\[
\tilde{t}_i = 
\begin{cases}
\min\left\{T_i, {\left\lfloor\frac{t_s}{2}\right\rfloor}\right\}, 
& i = 1, \ldots, k; \\
\max\left\{t_{i}, {\left\lceil\frac{t_s}{2}\right\rceil}\right\},
& i = k+1, \ldots, s.
\end{cases}
\]
 Then
\[
e_1(M) \leq \frac{\beta_0(M)}{(s+1)!} \prod_{i=1}^s \tilde {t}_i
\sum_{i=1}^s(\tilde {t}_i-i).
\]
\end{theorem}
When the resolution is quasi-pure, $t_i\ge T_{i-1}$ and this coincides with
the bounds in \cite{ES}.  
 Other authors gave also bounds for the first Hilbert coefficient $e_1$, see \cite{HH}, \cite{RV1}, \cite{RV2}, and \cite{E2}. We compare our bounds to that of Rossi-Valla in \cite{RV1} \cite{RV2}, and Elias in \cite{E2}, then conjecture that our result can be extended to all coefficients. 

\subsection*{Acknowledgements}
The authors also thank the University of Michigan for its hospitality where part of this work was done when the authors were visiting for the conference in honor of Craig Huneke. The second author is partially supported by a grant from the Infosys Foundation and a travel grant from the National Board of Higher Mathematics, India.  We thank the referee for a thorough reading of the manuscript and helpful comments.

\section{Notation}
Let $\Bbbk$ be a field and $R = \Bbbk[x_1, \ldots, x_n]$ be 
dimensional polynomial ring over $\Bbbk$ with $\deg x_i = 1$ for all 
$1 \leq i \leq n$.  Let $M$ be a finitely generated graded $R$-module.
Suppose $M$ is Gorenstein then the duality of the resolution can be
recorded as follow

If $\codim M=2k+1$ then\\

$0 \rightarrow  R(-d_s) \rightarrow \oplus_{j=1}^{\beta_1} R(-(d_s-d_{1j})
\rightarrow \ldots  \rightarrow \oplus_{j=1}^{\beta_k}R(-(d_s-d_{kj}))\\
\mbox{\hspace{7cm}} \rightarrow \oplus_{j=1}^{\beta_k}R(-d_{kj})
\rightarrow  \ldots \rightarrow \oplus_{j=1}^{\beta_1}R(-d_{1j})
\rightarrow R^{\beta_0}$
\\

and if codim $M=2k$ then\\
\\
$0\rightarrow  R(-d_s) {\rightarrow} \oplus_{j=1}^{\beta_1}R(-(d_s-d_{1j})
\rightarrow \ldots \rightarrow \oplus_{j=1}^{\beta_k/2=r_k}R(-(d_s-d_{kj}))
\oplus \oplus_{j=1}^{\beta_k/2=r_k}R(-d_{kj}) \\
 \mbox{\hspace {10cm}} \rightarrow  \ldots \rightarrow
 \oplus_{j=1}^{\beta_1}R(-d_{1j}) \rightarrow R^{\beta_0}$
 
 The minimal shifts in the resolution are:
$$t_i=\begin{cases}  \min_j d_{ij}& 1\leq i \leq k \\
d_s - \max_j d_{s-i,j}& k+1 \leq i < s \\
  d_s \hspace{0.25cm}& i = s 
\end{cases}$$

and the maximal shifts in the resolution are:
$$T_i=\begin{cases}
    \max_j  d_{ij}&  1\leq i \leq k \\
d_s - \min_jd_{s-i,j}&  k+1 \leq i < s \\
 d_s&  i = s 
\end{cases}$$

The \define{graded Betti numbers} of $M$ denoted by 
$\beta_{i,j}(M)$ are the number of copies of $R(-d_{ij})$ that appear at
homological degree $i$, in a minimal $R$-free resolution of $M$. We have 
$$\sum_j\beta_{i,j}= \begin{cases} \beta_i &i \leq k\\
\beta_{s-i}& k+1 \leq i < s \\
1& i=s
\end{cases}$$ 
We think of the collection
$\{\beta_{i,j}(M) : 0 \leq i \leq n, j \in \ints\}$ as an element
\[
\beta(M) = \left(\beta_{i,j}(M)\right)_{0 \leq i \leq n, j \in \naturals} 
\in \bbB := 
\bigoplus\limits_{i=0}^n \bigoplus\limits_{j \in  \mathbb{Z}} \mathbb{Q},
\]
and call it the \emph{Betti table} of $M$.  
In general, 
a \define{rational Betti table} $\beta$ is an element
$\beta = \left(\beta_{i,j}\right)_{0 \leq i \leq n, j \in \naturals} \in
\bbB$ 
such that: 
\begin{asparaenum}
\item
for all $0 \leq i \leq n$, $\beta_{i,j} = 0$ for finitely many $j$,
\item for all $i>0$ and for all $j$, if 
$\beta_{i,j} \neq 0$ then there exists $j' < j$ such that $\beta_{i-1,j'}
\neq 0$.
\end{asparaenum}
Let
$\beta = \left(\beta_{i,j}\right)_{0 \leq i \leq n, j \in \naturals}$ be a
rational Betti table.
Its \define{length} is $\max \{i : \beta_{i,j} \neq 0 \;\text{for some}\;
j\}$.

\subsection*{Vandermonde matrices}

 Given $(\alpha_1, \alpha_2, \ldots , \alpha_k)$ a sequence of real
 numbers, we denote the following Vandermonde determinants by
$$V_t=V_t(\alpha_1, \alpha_2, \ldots  , \alpha_{k}) = \left|
\begin{array}{cccc}
  1 & 1& .... & 1 \\
  \alpha_1&\alpha_2& ...&\alpha_k \\
  \alpha_1^2& \alpha_2^2 & ... &\alpha_k^2 \\
  \vdots & \vdots & ...&\vdots  \\
  \alpha_1^{k-2} & \alpha_2^{k-2} & ... & \alpha_k^{k-2} \\
  \alpha_1^{k-1+t}&\alpha_2^{k-1+t}& ... &\alpha_k^{k-1+t}\\
\end{array}\right|$$
 
 $$= \displaystyle \prod_{1 \leq j<i \leq
 k}(\alpha_i-\alpha_j)\displaystyle \sum_{\gamma_1+\gamma_2 + \ldots
 + \gamma_k=t}(\alpha_1^{\gamma_1}.\alpha_2^{\gamma_2} \ldots
 \alpha_k^{\gamma_k})$$
 
We denote $V_0$ by $V$.
  \begin{remark} $V_t(\alpha_1, \alpha_2, \ldots  , \alpha_{k}) \geq 0 $ if
  the sequence is in ascending order.
\end{remark} 

A \define{degree sequence} of \define{length} $s$ is an increasing sequence
$d = (d_0 < d_1 < \cdots < d_s)$ of integers.  An $R$-module $M$ is said
to have a \define{pure resolution} of type $d$ if its Betti table 
$\beta_{ij} \neq 0 $ if and only if $j= d_i, 1\le i\le s$. 
By the Herzog-K\"uhl equations~\cite{HK},
$\beta(M)$ is a positive rational multiple of the pure Betti
table, which we denote by $\beta(d)$, given by:

 \begin{equation}\label{equation:HerzogKuhl}
\beta(d)_{i,j} = 
\begin{cases}
\frac{1}{\prod_{l\neq i} |d_l-d_i|}, & 0 \leq i \leq s \;\text{and}\; j=d_i
\\
0, & \;\text{otherwise}.
\end{cases}
\end{equation}

 We call the Betti
table $\beta(d)$ defined in~\eqref{equation:HerzogKuhl}, the \define{pure
Betti table} associated to $d$. For $0 \leq i \leq s$, 
write $\beta_i(d) = \beta(d)_{i,d_i}$.

Eisenbud, {Fl{\o}ystad} and Weyman~\cite[Theorem~0.1]{EFWpureResln07} (in
characteristic zero) and Eisenbud and
Schreyer~\cite[Theorem~0.1]{ES} showed that for all degree sequences $d$,
there is a Cohen-Macaulay $R$-module $M$ such that $\beta(M)$ is a rational
multiple of $\beta(d)$. Moreover, for all $R$-modules $M$,
$\beta(M)$ can be written as a non-negative rational combination of the
$\beta(d)$~\cite{BS08}; if we take a saturated chain of degree sequences,
(in the set of all degree sequences, a saturated chain with respect to the partially order given by by point-wise comparison)
then the
non-negative rational coefficients in the decomposition are
unique~\cite[Theorem~0.2]{ES}.

\subsection{Self-dual resolutions and symmetrized Betti tables}

Let $\beta$ be a Betti table.  Let $s$ and $N$ be integers.
We say that $\beta$ is 
\define{$(s,N)$-self-dual}   if $\beta_{i,j} = \beta_{s-i, N-j}$ for all
$i,j$. 
We say that $\beta$ is \define{self-dual} 
if there exist $s$ and $N$ such that $\beta$ is $(s,N)$-self-dual.
If $\beta$ is self-dual, then $s$ is the length of $\beta$ and
$N = \max\{j : \beta_{s,j} \neq 0\} + \min\{j : \beta_{0,j} \neq 0\}$.

\begin{definition} 
\label{definition:symmetrized}
Let  $d = (d_0 < \cdots < d_s)$ be a degree sequence and $N \geq d_0+d_s$.
Let $d^{\vee, N} = (N-d_s< \cdots < N-d_0)$. 
Denote the pure Betti table
associated to $d^{\vee, N}$ by $\beta^{\vee, N}(d)$. Similarly, set 
$\beta^{\vee, N}_i(d) = \beta^{\vee, N}(d)_{i, N-d_{s-i}}$. 
Let $\beta_{\mathrm{sym}}(d,N) = \beta(d) + \beta^{\vee, N}(d)$. 
We call $\beta_{\mathrm{sym}}(d,N) $ the \define{symmetrized pure Betti
table}, given by symmetrizing $d$ with respect to $N$.
\end{definition}

Peskine and Szpiro~\cite{PS} showed that for a finitely generated
graded $R$-module $M$ with $s=\codim M$,
\[
\sum_{i=0}^{\projdim M} (-1)^i\displaystyle \sum_{j} d_{ij}^t = 
\begin{cases}
0 & \;\text{if}\; 0 \leq t < s \\  
(-1)^ss!e(M) & \;\text{if}\; t = s.
\end{cases}
\]

Moreover a similar expression to all coefficients was given in \cite[Lemma
3.2]{E} and \cite[Theorem 3.2]{ES1}, for every finitely generated
graded $R$-module $M$ of codimension $s$. We have,
$$
\displaystyle \sum_{r=0}^{l}(-1)^{l-r}\nu_{l-r} \displaystyle\sum
_{i=0}^{\projdim M}(-1)^i  \displaystyle\sum_{j=1}^{b_i} d_{ij}^{t+r} = 
(-1)^s(s+l)!e_l(M)  \;\text{if}\; t = s.
$$
with $\nu_{l-r}=\displaystyle \sum_{\tiny{1 \leq \xi_1 <\xi_2 <\cdots
<\xi_{l-r} \leq s+l-1}}\xi_1 \xi_2 \cdots  \xi_{l-r}$ and $\nu_0=1.$

Let $d$ be a degree sequence of length $s$.
Since the pure Betti table $\beta(d)$ is, up to
multiplication by a rational number, the Betti table of a Cohen-Macaulay
$R$-module of codimension $s$, we see  that 
$\sum_{i=0}^s (-1)^i\beta_{i}(d)d_i^l = 0$ for all $0 \leq l < s$.
Further, by a direct computation, we can see that $\sum_{i=0}^s
(-1)^i\beta_{i}(d)d_i^{s+r} = (-1)^s\frac{V_r(d_1,\ldots d_s)}{V(d_1,\ldots
d_s)}= (-1)^s \displaystyle \sum_{ \tiny {\alpha_1+\ldots
\alpha_s=r}}~d_1^{\alpha_1}d_2^{\alpha_1}\ldots d_s^{\alpha_s}$. 
\begin{definition} Let $d=(d_1, \ldots, d_s) $ be a degree sequence and $\beta (d)$ and $\beta(d^{\vee, N})$ be as in Definition \ref {definition:symmetrized}.  Then we define, $e_l$ as follows:
$$(s+l)!e_l(\beta(d)) =  \displaystyle
\sum_{r=0}^{l}(-1)^{l-r}\nu_{l-r}\displaystyle \sum_{\alpha_1+\ldots
\alpha_s=r}~d_1^{\alpha_1}d_2^{\alpha_2}\ldots d_s^{\alpha_s} $$
and 
$$(s+l)!e_l (\beta(d^{\vee, N})) = \displaystyle
\sum_{r=0}^{l}(-1)^{l-r}\nu_{l-r}\displaystyle \sum_{\tiny{\alpha_1+\ldots
\alpha_s=r}}(N- d_{s-1})^{\alpha_1}(N-d_{s-2})^{\alpha_2}\ldots
(N-d_0)^{\alpha_{s}}$$
\end{definition}
\vskip .2truein
 Note that the above sums can be factored as  follows
\begin{lemma} For any sequence of integers $y_1 < y_2< \ldots <y_s$
$$\displaystyle \sum_{r=0}^{l}(-1)^{l-r}\nu_{l-r}\displaystyle \sum_{\tiny
{\alpha_1+\ldots \alpha_s=r}}~y_1^{\alpha_1}y_2^{\alpha_2}\ldots
y_s^{\alpha_s} =\displaystyle \sum_{1 \leq i_1 \leq \cdots  i_l \leq
s}\displaystyle \prod_{t=1}^l (y_{i_t}-(i_{t}+t-1)) $$
where $  \nu_{l-r}=\displaystyle \prod_{1 \leq \beta_1 < \cdots
<\beta_{l-r} \leq s+l-1}\beta_1\cdots \beta_{l-r} $.
 \end{lemma}
This has been proved in ~\cite[Lemma~4.8]{ES1} and   \cite [Theorem
4.2]{E}. We repeat it for the sake of completeness.

\begin{proof}
For any given tuples $1  \leq \gamma_1 \leq \cdots \leq \gamma_r \leq s$
with $0 \leq r \leq l$, we have $1 \leq \beta_1 \leq \cdots \leq
\beta_{l-r} \leq s$ such that $\left\{ \gamma_1, \cdots \gamma_r \right\}
\cup \left\{ \beta_1, \cdots \beta_{l-r} \right\}=\left\{i_1, \cdots i_l
\right\}$.  In the product $\displaystyle \sum_{1 \leq i_1 \leq \cdots  i_l
\leq s}\displaystyle \prod_{t=1}^l (y_{i_t}-(i_{t}+t-1))$, the coefficient
of $\displaystyle \prod_{1\leq  \gamma_1 \leq \cdots \leq \gamma_l \leq
s}y_{\gamma_t}$  is $ \displaystyle \prod_{1 \leq \beta_1 < \cdots
<\beta_{l-r} \leq s+l-1}\beta_1\cdots \beta_{l-r}  = \nu_{l-r}$ since
$i_t+t-1$ is strictly increasing until $i_l+l-1$.
\end{proof}

Since 
$\beta_{\mathrm{sym}}(d,N) = \beta(d) + \beta^{\vee, N}(d)$, we see that 
\begin{equation} \label{equation:multPureBettiTable}(s+l)!e_l
(\beta_{\mathrm{sym}}(d,N))= \displaystyle \sum_{1 \leq i_1 \leq \cdots
i_l \leq s}\displaystyle \prod_{t=1}^l (d_{i_t}-(i_{t}+t-1))+\displaystyle
\sum_{1 \leq i_1 \leq \cdots  i_l \leq s}\displaystyle \prod_{t=1}^l
(N-d_{s-i_t}-(i_{t}+t-1))\end{equation} 

The Betti table of a Gorenstein module can be
decomposed into a non-negative rational combination of symmetrized pure
Betti tables, by   \cite{EKS}

\begin{proposition}\cite {EKS}[Proposition 2.4]
\label{proposition:moduleBettiDec}
Let $M$ be finitely generated graded Cohen-Macaulay $R$-module with $\codim
M = s$, generated minimally by homogeneous elements of degree zero. 
Suppose that $\beta(M)$ is self-dual. Let $N=T_s=t_s$.
Then there exist degree sequences $d^\alpha, 0 \leq \alpha \leq a$ for some
$a \in \naturals$ and positive rational numbers $r_\alpha, 0 \leq \alpha
\leq a$ such
that 
\[
\beta(M) = \sum_{\alpha=0}^a r_\alpha\beta_{\mathrm{sym}}(d^\alpha, N).
\]
Moreover,
\begin{enumerate}  
\item the $d^\alpha$ are degree sequences of length $s$ and they are not
$(s,N)$-dual to each other.
\item $d^{\alpha+1} > d^{\alpha}$ for all $0 \leq \alpha \leq a-1$.
\item $N \ge d^{\alpha}_i+d^{\alpha}_{s-i} $ for all $\alpha$ and $i$, or
equivalently, $d^{\alpha} \leq (d^{\alpha})^{\vee, N}$ for all $\alpha$.
\end{enumerate}
\end{proposition}

As a consequence, there is a formula for all the Hilbert coefficients,
$e_{\ell}$.   Let $M, N, d^{\alpha}$ be as in the theorem above. 

\[
e_{\ell} (M) = \sum_{\alpha=0}^a r_\alpha\beta_{\mathrm{sym}}(d^\alpha, N).
\]

Using this, one can get a stronger upper bound for the multiplicity $e_0$
of Gorenstein modules \cite{EKS}[Theorem 3.1].  
That is, 
   \[
e(M) \leq \frac{\beta_0(M)}{s!}
\prod_{i=1}^k \min\left\{T_i,
{\left\lfloor\frac{t_s}{2}\right\rfloor}\right\}
\prod_{i=k+1}^s \max\left\{t_{i}, 
{\left\lceil\frac{t_s}{2}\right\rceil}\right\}.
\]

Again, when the resolution is quasi-pure this generalizes the bound in
\cite{H}.   
This paper is an attempt to generalize this to higher Hilbert coefficients.
We prove an analogous bound for the first coefficient, $e_1$ and conjecture
a bound for the  higher coefficients.$e_i, i\ge 2$.

\section{Upper bound for  $e_1$}
\label{sec:proof}
 
\begin{notation} 
We define some notation used in the statements and proofs of our results.
Let $d  = (0, d_1, \ldots, d_s)$ be a non-decreasing sequence of
integers.
\begin{compactenum}

\item For a sequence 
$d' = (0, d'_1, \ldots, d'_s)$ of integers, we write $d < d'$ if $d \neq
d'$ and $d_i \leq d'_i$ for every $1 \leq i \leq s$.

\item $\tilde{d} = (\tilde{d}_1, \ldots, \tilde d_s)$ be the
(non-decreasing) sequence
\[
\tilde{d}_i = 
\begin{cases}
\min\left\{d_s-d_{s-i}, {\left\lfloor\frac{d_s}{2}\right\rfloor}\right\}, 
& i = 1, \ldots, k; \\
\max\left\{d_{i}, {\left\lceil\frac{d_s}{2}\right\rceil}\right\},
& i = k+1, \ldots s.
\end{cases}
\]

\item Suppose that $d_s \geq d_i + d_{s-i}$ for every $0 \leq i \leq s$.
Write
\begin{compactenum}
\item $b_{d} = \beta_0(d) + \beta_0(d^{\vee,d_s})$.
\item $\Psi_{d} = \prod_{i=1}^s \tilde d_i$.
\end{compactenum}

\item For $1 \leq l \leq s$, define
\[
f_l(\tilde d)=f_l(\tilde d_1,  \ldots,\tilde d_s)=\sum_{1 \leq i_1 \leq
\cdots \leq i_l \leq s}
\prod_{t=1}^l (\tilde d_{i_t}-(i_{t}+t-1)).
\]
Set $f_0(\tilde d_1,  \ldots, \tilde d_s) =1$.
\end{compactenum}

\end{notation}

This section is devoted to the proof of Theorem \ref{upperboundtheorem},
which is about finding 
an
upper bound for the first Hilbert
coefficient of Gorenstein algebras. 

\begin{theorem}\label{upperboundtheorem} Let $M$ be a finitely generated
graded Gorenstein
$R$-module, minimally generated by homogeneous elements of degree zero.
Let $s = \codim M$ and $k = \lfloor \frac{s}{2} \rfloor$.
Let $\beta_0(M)$ denote the minimal number of generators of $M$.
For $0 \leq i \leq s$, write 
$t_i = t_i(M) = \min\{j : \Tor_i^R(\Bbbk, M)_j \neq 0\}$ and $m =
(t_i)_{i=0, \ldots, s}$. Then
\[
e_1(M) \leq \frac{\beta_0(M)}{(s+1)!} \Psi_m f_1(\tilde t).
\]
\end{theorem}

This upper bound coincides with the upper bound of Gorenstein algebras with
quasi-pure resolutions found in~\cite{ES1}.
In order to prove this theorem, we find the upper bound for $e_
1(\beta_{\mathrm{sym}}(d,N))$ of a symmetrized pure Betti table then use
the Boij and S\"oderberg decomposition in order to generalize our result to
$e_1(M)$. To proceed, we first need the following lemma,

\begin{lemma}
\label{lemma:comparisonBettiPsi}
Let $d$ and $d'$ be degree sequences such that $d_0 = 0$ and
$d < d' \leq (d')^{\vee, d_s} < d^{\vee, d_s}$.
Then $\Psi_{d} f_1(\tilde{d})\geq \Psi_{d'}f_1(\tilde{d'})$.
\end{lemma}

\begin{proof}By induction on $\sum_{i} d'_i-d_i$,  we may assume,  
without loss of generality, that there exists $j$ such that
$d'_j = d_j+1$ and $d'_i = d_i$ for all $i \neq j$.
Moreover, if $1 \leq  i \leq s-k-1$, then $d_i$
does not figure in the expression for $\Psi_{d} f_l(\tilde{d})$, so
we may assume that $j \geq s-k$. Additionally, $j \leq s-1$.
We rewrite $\Psi_d$ as
\begin{equation}
\label{equation:PsidReWritten}
\Psi_d = 
\prod_{i=s-k}^{s-1} \min\left\{d_s-d_i, 
{\left\lfloor\frac{d_s}{2}\right\rfloor}\right\}
\prod_{i=k+1}^s \max\left\{d_{i}, 
{\left\lceil\frac{d_s}{2}\right\rceil}\right\}.
\end{equation}
We note that 
\begin{align*}
f_1(\tilde{d}) & = \sum_{1 \leq i_1  \leq s} (\tilde{d_{i_1}}-i_{1}) =
\sum_{1 \leq i  \leq s}\tilde{ d_{i}} - \binom{s+1}{2}
\end{align*}

Two cases arise: $j < k+1$ and $j \geq k+1$. The first case is possible if
and only if $s=2k$ and $j=k$.
In this case, $d_k$ appears only once
in~\eqref{equation:PsidReWritten}, and since $d_k \leq
d_{s}-d_{s-k} = d_s-d_k$, we get 
$d_s - d_k \geq {\left\lfloor\frac{d_s}{2}\right\rfloor}$.
By the hypothesis that $d' \leq (d')^{\vee, d_s}$, 
$d_s - d_k - 1 \geq d_k+1$, so 
$d_s - d_k -1 \geq {\left\lfloor\frac{d_s}{2}\right\rfloor}$. We get
$\min\left\{d_s-d_k,
{\left\lfloor\frac{d_s}{2}\right\rfloor}\right\}=\min\left\{d_s-d_k-1,
{\left\lfloor\frac{d_s}{2}\right\rfloor}\right\}=
{\left\lfloor\frac{d_s}{2}\right\rfloor}$. In the expression
$f_1(\tilde{d'})$, only $d_s-d_k -1 $ is involved, so we get
$f_1(\tilde{d'})\Psi_ {d'}= f_1(\tilde{d})(\Psi_{d})$.

In the second case (i.e., $j \geq k+1$), $d_j$ appears twice in 
in~\eqref{equation:PsidReWritten}. We need to show that 
\[
\frac{\Psi_{d'}f_l(\tilde{d'})}{\Psi_d f_l(\tilde{d})} \leq 1 \]

If $d_j < \left\lceil\frac{d_s}{2}\right\rceil$, then $d_s - d_j >
\left\lfloor\frac{d_s}{2}\right\rfloor$. So we get that $d_j+1 \leq
\left\lceil\frac{d_s}{2}\right\rceil$ and $d_s - d_j-1 \geq
\left\lfloor\frac{d_s}{2}\right\rfloor$. Hence, $\min\left\{d_s-d_j,
{\left\lfloor\frac{d_s}{2}\right\rfloor}\right\}=\min\left\{d_s-d_j-1,
{\left\lfloor\frac{d_s}{2}\right\rfloor}\right\}=
{\left\lfloor\frac{d_s}{2}\right\rfloor}$ and $\max\left\{d_{j},
{\left\lceil\frac{d_s}{2}\right\rceil}\right\}=\max\left\{d_{j},
{\left\lceil\frac{d_s}{2}\right\rceil}\right\}={\left\lceil\frac{d_s}{2}\right\rceil}$.
As a result we obtain $f_1(\tilde{d'})\Psi_ {d'}=
f_1(\tilde{d})(\Psi_{d})$.

Let us now consider the case where $d_j \geq
\left\lceil\frac{d_s}{2}\right\rceil$ so $d_j+1 \geq
\left\lceil\frac{d_s}{2}\right\rceil $. This implies that $d_s-d_j-1 <
d_s-d_j \leq  \left\lfloor\frac{d_s}{2}\right\rfloor$. It suffices to show
that \[ \frac{\Psi_{d'}f_l(\tilde{d'})}{\Psi_df_l(\tilde{d})} = 
\frac{(d_s - d_j - 1)(d_j+1)f_1(\tilde{d'})}{(d_s - d_j)d_jf_1(\tilde{d})}
\leq 1. \]

We compute the difference between the numerator and the denominator. We
write

 \begin{equation} \label{equation1}
(d_s - d_j - 1)(d_j+1)f_1(\tilde{d'})-(d_s -
d_j)d_jf_1(\tilde{d})=-(2d_j-d_s+1)f_1(\tilde{d})+
(d_s-d_j-1)(d_j+1)(f_1(\tilde{d'})-f_1(\tilde{d})) \end{equation}
 where 
 $$ f_1(\tilde{d}) = \sum_{1 \leq i  \leq s}\tilde{ d_{i}} -
 \binom{s+1}{2}=s \left\lfloor\frac{d_s}{2}\right\rfloor - \binom{s+1}{2}$$
and
$$f_1(\tilde{d'}) =  \sum_{1 \leq i  \leq s}\tilde{ d'_{i}} -
\binom{s+1}{2}=s \left\lfloor\frac{d_s}{2}\right\rfloor -
2-\binom{s+1}{2}$$

So Equation \ref{equation1} is equal to:
$$-(2d_j-d_s+1)(s \left\lfloor\frac{d_s}{2}\right\rfloor - \binom{s+1}{2})+
(d_s-d_j-1)(d_j+1)(-2)$$

We know that $d_j \geq \left\lceil\frac{d_s}{2}\right\rceil$ and $s
\left\lfloor\frac{d_s}{2}\right\rfloor -\frac{s(s+1)}{2}=s(
\left\lfloor\frac{d_s}{2}\right\rfloor-\frac{s+1}{2})\geq 0$ since $d_s
\geq s+2$.
Hence $(d_s - d_j - 1)(d_j+1)f_1(\tilde{d'})-(d_s -
d_j)d_jf_1(\tilde{d})\leq 0$, and the proof is done. 
\end{proof}

The next proposition is needed for the proof of
Theorem~\ref{upperboundtheorem}.
Let $d  = (0, d_1, \ldots, d_s)$ be a degree sequence such that $d_s \geq 
d_i + d_{s-i}$ for all $0 \leq i \leq s$. We show that the first
coefficient $e_1 (\beta_{\mathrm{sym}}(d,d_s)$ of symmetrized pure
sequences satisfy the upper bound in our main theorem. When $l=1$, we write
from equation \ref {equation:multPureBettiTable}:

\begin{equation}\label{equ2}
  \begin{array}{ccc}(s+l)!e_1 (\beta_{\mathrm{sym}}(d,d_s))= \displaystyle
  \sum_{1 \leq i_1 \leq s}\displaystyle  
(d_{i_1}-i_{1})+\displaystyle \sum_{1 \leq i_1 \leq
s}\displaystyle  (d_s-d_{s-i_1}-i_{1})  \\

=f_1(d)+f_1(d^{\vee, d_s})
\end{array}
\end{equation}

We note that 
\begin{align*}
f_1(d) &= \sum_{1 \leq i  \leq s} d_{i} - \binom{s+1}{2}
; \\
f_1(d^{\vee, d_s}) & = sd_s - \sum_{1 \leq i  \leq s} d_{i} -
\binom{s+1}{2}.
\end{align*}

\begin{proposition}
\label{proposition:theoremforsymmtrizedpure}

Let $d  = (0, d_1, \ldots, d_s)$ be a degree sequence such that 
$d \leq d^{\vee, d_s}$. Then 
\[
f_1(d) + f_1(d^{\vee, d_s}) \leq b_d \Psi_{d} f_1(\tilde{d}).
\]

\end{proposition}

\begin{proof} 
We again prove this by induction on 
$\sum_{i} \left((d^{\vee, d_s})_i - d_i\right)=
\sum_{i} (d_s-d_{s-i}-d_i)$, which is non-negative by our hypothesis.
If $\sum_{i} (d_s-d_{s-i}-d_i) = 0$ (equivalently, $d = d^{\vee, d_s}$), 
then $\tilde{d} = d$, so the assertion follows from noting that 
$b_{d} \Psi_{d}  \geq 2$~\cite[Proposition~3.5]{EKS}.

If $d < d^{\vee, d_s}$, then there exists $j \geq k$ such that 
$d_j < d_{s}-d_{s-j}$. Pick $j$ to be maximal with this property. 
Note that $j < s$, so $d'_s = d_s$. 
Since $d_j < d_s - d_{s-j}$, $d_{j+1} = d_s - d_{s-j-1}$ and 
$d_{s-j} > d_{s-j-1}$, we see that
$d' := (0, d_1, \cdots, d_{j-1}, d_j+1, d_{j+1}, \cdots, d_s)$ is a degree
sequence, that
$d' \leq (d')^{\vee,d_s}$ 
and that
$\sum_{i} (d_s-d_{s-i}-d_i) > \sum_{i} (d_s-d'_{s-i}-d'_i)$. Hence by 
induction
\[
f_1(d') + f_1((d')^{\vee, d_s}) \leq b_{d'} \Psi_{d'} f_1(\tilde{d'}).
\]

Therefore $f_1(d') + f_1((d')^{\vee, d_s}) = 
f_1(d) + f_1(d^{\vee, d_s})$. We now show that 
$f_1(\tilde{d'}) = f_1(\tilde{d})$. To this end, we consider various cases:

\begin{enumerate}

\item
\underline{$j < s-k$}: Then $\tilde{d'} = \tilde{d}$, so 
$f_1(\tilde{d'}) = f_1(\tilde{d})$. 

\item
\underline{$j = s-k$}: We consider two sub-cases. 
\begin{enumerate}
\item
\underline{$s = 2k$}: Then 
\[
(\tilde{d'})_i = 
\begin{cases}
\tilde{d}_i, & 1 \leq i \leq s, i \neq k \\
\min\left\{d_s-d_{k}-1, {\left\lfloor\frac{d_s}{2}\right\rfloor}\right\}, 
& i = k \\
\end{cases}
\]
Note that since $d_k < d_s - d_k$, 
$d_s - d_k-1 \geq {\left\lfloor\frac{d_s}{2}\right\rfloor}$, so 
$(\tilde{d'})_k = \tilde{d}_k = {\left\lfloor\frac{d_s}{2}\right\rfloor}$.
Therefore $d' = d$ and hence 
$f_1(\tilde{d'}) = f_1(\tilde{d})$. 

\item
\underline{$s = 2k+1$}: Then 
\[
(\tilde{d'})_i = 
\begin{cases}
\tilde{d}_i, & 1 \leq i \leq s, i \neq k \;\text{and}\; i \neq k+1 \\
\min\left\{d_s-d_{k+1}-1, {\left\lfloor\frac{d_s}{2}\right\rfloor}\right\},
& i = k \\
\max\left\{d_{k+1}+1, {\left\lceil\frac{d_s}{2}\right\rceil}\right\},
& i = k+1 \\
\end{cases}
\]
Therefore
$f_1(\tilde{d'}) = f_1(\tilde{d})$. 
\end{enumerate}

\item
\underline{$j > s-k$}: Then 
\[
(\tilde{d'})_i = 
\begin{cases}
\tilde{d}_i, & 1 \leq i \leq s, i \neq j \;\text{and}\; i \neq s-j \\
\min\left\{d_s-d_{j}-1, {\left\lfloor\frac{d_s}{2}\right\rfloor}\right\},
& i = s-j \\
\max\left\{d_{j}+1, {\left\lceil\frac{d_s}{2}\right\rceil}\right\},
& i = j \\
\end{cases}
\]
Therefore
$f_1(\tilde{d'}) = f_1(\tilde{d})$. 
\end{enumerate}

Now note that $b_d \Psi_d \geq b_{d'} \Psi_{d'}$~\cite[Proof of
Proposition~3.5, p.126]{EKS}. This completes the proof of the Proposition.
\end{proof}

\begin{proof}[Proof of Theorem~\ref{upperboundtheorem}]
Pick degree sequences $d^\alpha$ and non-negative rational numbers 
$r_\alpha$ as in Proposition~\ref{proposition:moduleBettiDec}. 
We need to show that $(s+1)!e_1(M) \leq \beta_0(M)\Psi_{t}f_1(\tilde{t})$.
We get this as follows:

\begin{equation}
\begin{split}
(s+1)!e_1(M)
& = (s+1)! \sum_{\alpha} r_{\alpha} e_1(\beta_{sym}(d^{\alpha})) \\
& = \sum_{\alpha} r_{\alpha} \left(f_1(d^\alpha) +
f_1((d^\alpha)^{\vee, d_s})\right) 
\;\; (\text{by}~\eqref{equ2}\\
& \leq \sum_{\alpha} r_{\alpha} b_{d^\alpha} \Psi_{d^\alpha}
f_1(\tilde{d^\alpha}) 
\;\; (\text{by Proposition}~\ref{proposition:theoremforsymmtrizedpure})\\
& \leq \left(\sum_{\alpha} r_{\alpha} b_{d^\alpha} \right)
\Psi_{t} f_1(\tilde{t})
\;\; (\text{by Lemma}~\ref{lemma:comparisonBettiPsi})\\
& = \beta_0(M) \Psi_{t} f_1(\tilde{t}).
\end{split}
\end{equation}
\end{proof}

Next we  give examples for the bound that we found and we compare our
result to the bound found by Rossi-Valla in \cite[Theorem 3.2]{RV2} and Elias in \cite [Theorem 2.3]{E2}. We obtain a
much sharper bound. 

\begin{example} Let $R= \Bbbk[x,y,z,w,t,u,v]$ and $I=(y^2, z^2, w^2,
x^3-zwt, x^2y-wt^2, x^2z, xyz-x^2t, yzt-xt^2, yt^2, zt^2, t^3) $ a
homogeneous Gorenstein ideal of $R$. Note that  $e_1(R/I)= 90$, and the
Betti diagram is given as follow:
\tiny {\begin{verbatim} 
             0  1  2  3  4 5
      total: 1 11 28 28 11 1
          0: 1  .  .  .  . .
          1: .  3  .  .  . .
          2: .  8 20  8  . .
          3: .  .  8 20  8 .
          4: .  .  .  .  3 .
          5: .  .  .  .  . 1
\end{verbatim}}

\normalsize
We have $T_1 = 3, T_2=5, t_3=5 , t_4=7$, and $d_5=T_5=t_5=10$ and
$\frac{d_5}{2}=5$. By theorem \ref{upperboundtheorem},  $e_1(R/I) \leq
\frac{1}{6!}f_1(T_1, T_2, t_3, t_4, t_5) T_1 T_2t_3t_4t_5$, which implies
that $$e_1(R/I) \leq \frac{1}{6!}(2+3+2+3+5)3.5.5.7.10= 109.375$$  The bound found by  Elias in \cite[Theorem 2.3]{E2} is equal to  ${e_0 \choose 2} - {\mu(I)-d \choose 2}$ where $\mu(I)=11$ is the minimal number of generator of $I$ and $d=7$ the dimension of $R$. Since $e_0=90$, then the bound in \cite{E2} gives $e_1(R/I) \leq 624$.

\end{example}

In the above example, the minimal free resolution of $R/I$ is quasi-pure.
This case was done in \cite{ES1}. In the next example, $R/I$ has a  non
quasi-pure minimal free resolution.  

\begin{example}\label{ex3} Let $R= \Bbbk[x,y,z,w,s,t, u, v]$ and $I=(wt, zt, xt, zw, yw, xw, yz,
xz, x^2y+yt^2, x^3-y^3+yt^2-t^3, w^5+t^5, z^5-t^5)$ a homogeneous
Gorenstein ideal of $R$. Note that  $e_1(R/I)= 65$, and the Betti diagram
is given as follow:
 \tiny {\begin{verbatim} 
             0  1  2  3  4 5
      total: 1 12 29 29 12 1
          0: 1  .  .  .  . .
          1: .  8 14  9  2 .
          2: .  2  4  2  . .
          3: .  .  2  4  2 .
          4: .  2  9 14  8 .
          5: .  .  .  .  . 1
\end{verbatim}}
\normalsize
We have $T_1 = 5, T_2=6, t_3= 4, t_4=5$, and $d_5=T_5=t_5=10$ and
$\frac{d_5}{2}=5$. By theorem \ref{upperboundtheorem},  $e_1(R/I) \leq
\frac{1}{6!}f_1(T_1, \frac{d_5}{2}, \frac{d_5}{2}, t_4, t_5) T_1
\frac{d_5}{2} \frac{d_5}{2}t_4t_5$, which implies that $$e_1(R/I) \leq
\frac{1}{6!}(4+3+2+1+5)5.5.5.5.10 \simeq 130.20833$$ Whereas the bound
found by  Rossi-Valla in \cite[Theorem 3.2]{RV2} is equal to  ${e_0 \choose 2} - {\mu(I)-d \choose 2}-\lambda(R/I)+1 $ where $\mu(I)=12$ is the minimal number of generator of $I$, $d=8$ the dimension of $R$, and $\lambda(R/I)=e_0-e_1+e_2$. We note that $e_0(R/I)=26$ and $e_2(R/I)=68$. Hence the bound in \cite{RV2} gives $e_1(R/I) \leq 285$.
\end{example}

 Huneke and Hanumanthu \cite{HH} find a bound  for $e_1$ in a slightly different setting.  For a CM local ring $(R,m)$  of dim d and an ideal $I$  contained in $m^k$, for some $k>2$, let  $\lambda(R/I^{n+1}) = \sum _{i}(-1)^i e_i(I) {(n+d-i\choose d)}$.  In \cite{HH} corollary 3.7, they show that
$e_1 \le {((e_0-k)\choose 2)}$.  
\vskip .2truein
Finally, we conjecture that the above upper bound can be extended to all
Hilbert coefficients $e_i$'s. We believe the following is true.

\begin{conjecture} Let $M$ be a finitely generated graded Gorenstein
$R$-module, minimally generated by homogeneous elements of degree zero.
Let $s = \codim M$ and $k = \lfloor \frac{s}{2} \rfloor$.
Let $\beta_0(M)$ denote the minimal number of generators of $M$ and write
$m =( t_i)_{i=0, \ldots, s}$. Then
\[
e_j(M) \leq \frac{\beta_0(M)}{(s+1)!} \Psi_tf_j(\tilde t).
\]for all $0 \leq j \leq d-1$.
\end{conjecture}

\begin{example} In Example \ref {ex3}, we considered the ideal $I=(wt, zt, xt, zw, yw, xw, yz,
xz, x^2y+yt^2, x^3-y^3+yt^2-t^3, w^5+t^5, z^5-t^5)$ where $T_1 = 5, T_2=6, t_3=4 , t_4=5$, $d_5=T_5=t_5=10$ and
$\frac{d_5}{2}=5$. We have $e_2(R/I)=68$, and we note that 
 $$e_2(R/I) \leq \frac{1}{6!}150.5.5.5.5.10 \simeq 1302 $$ 

\end{example}
\providecommand{\bysame}{\leavevmode\hbox to3em{\hrulefill}\thinspace}
\providecommand{\MRhref}[2]{%
  \href{http://www.ams.org/mathscinet-getitem?mr=#1}{#2}
}
\providecommand{\href}[2]{#2}

\end{document}